\documentclass[12pt]{amsart}
\usepackage{amsmath,amssymb,amsfonts,amsthm,amstext,verbatim,url}
\usepackage{graphicx}
\RequirePackage[colorlinks,citecolor=blue,urlcolor=blue]{hyperref}

\theoremstyle{plain}

\newtheorem{theorem}{Theorem}
\newtheorem{definition}[theorem]{Definition}

\newtheorem{example}[theorem]{Example}

\newtheorem{remark}[theorem]{Remark}
\newtheorem{question}[theorem]{Question}
\numberwithin{equation}{section}
\numberwithin{theorem}{section}
\numberwithin{figure}{section}

\newcounter{mycount}\newcounter{mycount2}
\newenvironment{romlist}{\begin{list}{\rm(\roman{mycount})}%
   {\usecounter{mycount}\labelwidth=1cm\itemsep 0pt}}{\end{list}}

\newenvironment{letlist}{\begin{list}{\rm(\alph{mycount2})}%
   {\usecounter{mycount2}\labelwidth=1cm\itemsep 0pt}}{\end{list}}
   \newenvironment{Alist}{\begin{list}{\MakeUppercase{\alph{mycount}}.}%
   {\usecounter{mycount}\labelwidth=1cm\itemsep 0pt}}{\end{list}}

\newcommand\s{\sigma}

\newcommand\ol{\overline}
\newcommand\oo{\infty}
\newcommand\De{\Delta}
\newcommand\LL{{\mathbb L}}
\newcommand\HH{{\mathbb H}}
\newcommand\NN{{\mathbb N}}
\newcommand\TT{{\mathbb T}}
\newcommand\PP{{\mathbb P}}
\newcommand\EE{{\mathbb E}}

\newcommand\sP{{\mathcal P}}

\newcommand\sT{{\mathcal T}}
\newcommand\sA{{\mathcal A}}
\newcommand\sG{{\mathcal G}}

\newcommand\sN{{\mathcal N}}
\newcommand\sB{{\mathcal B}}

\newcommand\sQ{{\mathcal Q}}
\newcommand\ZZ{{\mathbb Z}}
\newcommand\RR{{\mathbb R}}
\newcommand\BB{{\mathbb B}}

\newcommand\wt{\widetilde}
\newcommand\om{\omega}

\renewcommand\a{\alpha}
\newcommand\Ga{\Gamma}

\newcommand\Si{\Sigma}
\newcommand\si{\sigma}
\newcommand\eps{\epsilon}

\newcommand\g{\gamma}
\newcommand\resp{respectively}
\newcommand\fish{F}

\newcommand\Stab{\text{\rm Stab}}

\newcommand\olG{{\ol G}}
\newcommand\olV{{\ol V}}
\newcommand\olE{{\ol E}}

\newcommand\vG{{\vec G}}

\newcommand\vE{{\vec E}}

\newcommand\vpi{\vec\pi}
\newcommand\vmu{\vec\mu}

\newcommand\q{\quad}
\newcommand\Aut{\text{\rm Aut}}
\newcommand\pd{\partial}

\newcommand\id{{\bold 1}} 

\newcommand\SLE{\text{\rm SLE}}
\renewcommand\O{\text{\rm O}}
\newcommand\muF{\mu^{\text{\rm F}}}
\newcommand\siF{\s^{\text{\rm F}}}
\newcommand\pc{p_{\text{\rm c}}}
\newcommand\Tc{T_{\text{\rm c}}}
\newcommand\ollG{G'}
\newcommand\ollE{E'}
\newcommand\TLF{\sT}
\newcommand\ghf{graph height function}
\newcommand\hdi{$\sH$-difference-invariant}
\newcommand\sH{{\mathcal H}}
\newcommand\qq{\qquad}
\newcommand\ughf{unimodular graph height function}

\newcommand\EGF{\mathrm{EFG}}
\newcommand\normal{\trianglelefteq}
\newcommand\EG{\mathrm{EG}}
\newcommand\mumin{\mu_{\mathrm{min}}}

\begin{document}
\title[Self-avoiding walks]{Self-avoiding walks and\\ connective constants}

\author[Grimmett]{Geoffrey R.\ Grimmett}
\address{Statistical Laboratory, Centre for
Mathematical Sciences, Cambridge University, Wilberforce Road,
Cambridge CB3 0WB, UK} 
\email{g.r.grimmett@statslab.cam.ac.uk}
\urladdr{\url{http://www.statslab.cam.ac.uk/~grg/}}
\address{Department of Mathematics, University of Connecticut, 341 Mansfield Road U1009, Storrs, CT 06269--1009, USA}

\email{zhongyang.li@uconn.edu}
\urladdr{\url{http://www.math.uconn.edu/~zhongyang/}}
\author[Li]{Zhongyang Li}

\begin{abstract}
The \emph{connective constant} $\mu(G)$
of a quasi-transitive graph $G$ is the asymptotic growth
rate of the number of self-avoiding walks (SAWs) on $G$
from a given starting vertex.
We survey several aspects of the relationship between the connective constant
and the underlying graph $G$. 
\begin{itemize}
\item We present upper and lower bounds for $\mu$ in terms
of the vertex-degree and girth of a transitive graph.
\item We discuss the question of whether $\mu\ge\phi$ for transitive cubic graphs
(where $\phi$ denotes the golden mean), and we introduce  the Fisher transformation for SAWs 
(that is, the replacement of vertices by triangles).
\item We present strict inequalities
for the connective constants $\mu(G)$ of transitive graphs $G$, as $G$ varies. 
\item As a consequence of the last, the connective constant of a Cayley graph 
of a finitely generated group decreases strictly when
a new relator is added, and increases strictly when a 
non-trivial group element is declared to be a further generator.
\item We describe  so-called \ghf s within an account of \lq bridges' for 
quasi-transitive graphs, and indicate that
the bridge constant equals the connective constant when the graph has a \ughf.
\item A partial answer is given to the question of the locality of connective constants, based
around the existence of \ughf s.
\item Examples are presented of Cayley graphs 
of finitely presented groups that possess  \ghf s (that are, in addition,
harmonic and unimodular), and that do not.
\item The review closes with a brief account of the \lq speed' of SAW.
\end{itemize}
\end{abstract}

\date{19 April 2017}

\maketitle

\vfill\eject
\tableofcontents

\vfill\eject
\section{Introduction}
\subsection{Self-avoiding walks}
\begin{figure}[htbp]
\centerline{\includegraphics[width=0.3\textwidth]{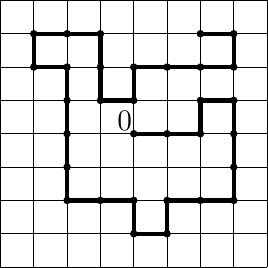}}
   \caption{A 31-step SAW from the origin of the square lattice.}
\label{fig:saws}
\end{figure}

A \emph{self-avoiding walk} (abbreviated to SAW) on a 
graph $G=(V,E)$ is a path that visits no vertex
more than once. An example of a  SAW on the square 
lattice is drawn in Figure \ref{fig:saws}.
SAWs were first introduced in the chemical theory of polymerization (see Orr \cite{Orr} and the book of Flory \cite{f}),
and their critical behaviour has attracted the 
abundant attention since of mathematicians and physicists 
(see, for example, the book of Madras and Slade \cite{ms}
and the lecture notes \cite{bdgs}). 

The theory of SAWs impinges on several areas of science including 
combinatorics, probability, and statistical mechanics.
Each of these areas poses its characteristic questions
concerning counting and geometry. The most fundamental problem is to count the
number of $n$-step SAWs from a given vertex, and this is the starting point
of a rich theory of geometry and phase transition.

Let $\s_n(v)$
be the number of $n$-step SAWs on $G$ starting  at the vertex $v$.
The following fundamental theorem of Hammersley asserts the existence of an
asymptotic growth rate for $\s_n(v)$ as $n \to \oo$. 
(See Section \ref{ssec:trans} for a definition of (quasi-)transitivity.)

\begin{theorem}\cite{jmhII}\label{jmh}
Let $G=(V,E)$ be an infinite, connected, quasi-transitive graph with finite vertex-degrees. There exists
$\mu=\mu(G)\in[1,\oo)$, called the \emph{connective constant} of $G$,  such that
\begin{equation}\label{connconst}
\lim_{n \to \oo} \s_n(v)^{1/n} = \mu, \qquad v \in V.
\end{equation}
\end{theorem}

At the heart of the proof is the observation by Hammersley and Morton \cite{hm} that 
(in the case of a transitive graph) $\log\sigma_n$ is a subadditive function. That is,
\begin{equation}\label{eq:subadditive}
\sigma_{m+n} \le \sigma_m \sigma_n, \qq m,n\ge 1.
\end{equation}

The value $\mu=\mu(G)$ depends evidently on the choice of graph $G$.
Indeed, $\mu$ may be viewed as  a \lq critical point', 
corresponding, in a sense, to the critical probability
of the percolation model, or the critical temperature of the Ising model.
Consider the generating function
\begin{equation}\label{eq:gf}
Z_v(x)=\sum_{w\in\Sigma(v)}x^{|w|}, \qq x \in \RR,
\end{equation}
where $\Sigma(v)$ is the set of finite SAWs starting from a given vertex $v$, and $|w|$ is the 
number of edges of $w$. 
Viewed as a power series, $Z_v(x)$ has radius of convergence $1/\mu$,
and thus a singularity at the point $x=1/\mu$.
Critical exponents may be introduced as in Section \ref{ssec:ce}.

In this paper we review certain properties of the connective constant $\mu(G)$, in particular
exact values (Section \ref{ssec:concon}), upper and lower bounds (Section \ref{sec:lb}), 
a sharp lower bound for cubic graphs, and the Fisher transformation
(Section \ref{sec:fisher}), strict inequalities
(Sections \ref{sec:si}--\ref{sec:cayley}), and the locality theorem (Section \ref{sec:loc}).
The results summarised here may be found largely in the work of the authors 
\cite{GrL2}--\cite{GrL7} and \cite{ZL16}. 
This review is an expanded and updated version of \cite{GL-revI}.

Previous work on SAWs tends to have been focussed on specific graphs such as 
the cubic lattices $\ZZ^d$ and certain two-dimensional lattices. 
In contrast, the results of \cite{GrL2}--\cite{GrL7} 
are directed at general classes of graphs that are \emph{quasi-transitive}, and often \emph{transitive}.
The work reviewed here may be the first systematic study of 
SAWs on general transitive and 
quasi-transitive graphs.
It is useful to have a reservoir of (quasi-)transitive graphs at one's disposal
for the construction and analysis of hypotheses, and to this end the Cayley
graphs of finitely generated groups play a significant role (see Section \ref{ssec:cayley}).
We note the recent result of Martineau \cite{Mart16}
that the set of connective constants of Cayley graphs contains
a Cantor space.

Notation for graphs and groups will be introduced when needed.
A number of questions are included in this review. 
The inclusion of a question does not of itself imply
either difficulty or importance. 

\subsection{Connective constants, exact values}\label{ssec:concon}
\begin{figure}[t]
 \centering
\includegraphics[width=0.8\textwidth]{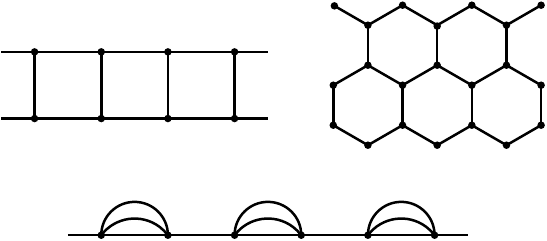}
  \caption{Three regular graphs: the ladder graph $\LL$,  
the hexagonal tiling $\HH$
of the plane, and the bridge graph $\BB_\Delta$ 
(with $\Delta=4$) obtained from $\ZZ$ by joining
every alternate pair of consecutive vertices by $\Delta-1$ parallel edges.}
  \label{fig:ladder-hex}
\end{figure}

For what graphs $G$ is $\mu(G)$ known exactly? 
There are a number of such graphs, 
which should be regarded as atypical in this regard.
We mention the \emph{ladder} $\LL$, the hexagonal lattice $\HH$, and the 
\emph{bridge graph}  $\BB_\De$ with degree $\De\ge 2$ of Figure \ref{fig:ladder-hex}, for which
\begin{equation}\label{2}
\mu(\LL) = \tfrac12(1 + \sqrt 5), \quad \mu(\HH) = \sqrt{2+\sqrt 2},
\quad \mu(\BB_\De) = \sqrt {\De-1}.
\end{equation}
See \cite[p.\ 184]{AJ90} and \cite{ds}
for the first two calculations. The third is elementary.

In contrast, the value of the
connective constant of the square grid $\ZZ^2$ is unknown, and a substantial
amount of work has been devoted to obtaining good bounds. 
The best rigorous bounds known currently to the authors are those of \cite{j04,PT}, namely
(to $5$ significant figures)
\begin{equation*}
2.6256\leq \mu(\ZZ^2)\leq 2.6792,
\end{equation*}
and more precise numerical estimates are available, including the estimate
 $\mu \approx 2.63815\dots$ of \cite{Jens03}.

We make some remarks about the three graphs of Figure \ref{fig:ladder-hex}.
There is a correspondence between the Fibonacci sequence 
and counts of SAWs on the ladder graph $\LL$ (see, for example \cite{Zeil}), whereby one obtains that 
$\mu(\LL)$ equals the golden ratio $\phi := \frac12(1+\sqrt 5)$.
We ask in Question \ref{qn:gmean} 
whether $\mu(G) \ge \phi$ for all infinite, simple, cubic, transitive graphs,
and we discuss evidence for a positive answer to this question.

Amongst a certain class of $\De$-regular 
graphs permitted to possess multiple edges, the bridge graph
$\BB_\De=\sqrt{\De-1}$ is extremal in the sense that $\mu(\BB_\De)$ is least.
See the discussion of Section \ref{sec:lb}.

The proof that $\mu(\HH)=\sqrt{2+\sqrt{2}}$ by 
Duminil-Copin and Smirnov \cite{ds}
is a very significant recent result. The value $\sqrt{2+\sqrt 2}$ 
emerged in the physics literature through work  of Nienhuis \cite{Nien}
motivated originally by 
renormalization group theory. Its proof in \cite{ds} is based
on the construction of an observable with certain properties of discrete holomorphicity,
complemented by a neat use of the bridge decomposition introduced by Hammersley and Welsh
\cite{HW62}. The bridge decomposition has been used since to prove the 
locality theorem for connective constants (see Section \ref{sec:bridge} 
and Theorem \ref{thm2}).

\subsection{Three problems on the square lattice}
There are a number of beautiful open problems associated with SAWs and
connective constants, of which we select three. Our first problem is to prove that
a random $n$-step SAW from the origin of $\ZZ^2$ converges, when suitably rescaled,
to the Schramm--Loewner curve $\SLE_{8/3}$.
This important conjecture has been discussed and formalized 
by Lawler, Schramm, and Werner \cite{lsw}. 

\begin{question}
Does a uniformly distributed  $n$-step SAW from the origin of $\ZZ^2$ converge,
when suitably rescaled,
to the random curve $\SLE_{8/3}$?
\end{question}

Recent progress in this direction was made by  Gwynne and Miller  \cite{GwM}, who 
proved that a  SAW on a random quadrangulation converges
to $\SLE_{8/3}$ on a certain Liouville-gravity surface.

There is an important class of results usually referred to as the `pattern theorem'.
In Kesten's original paper \cite{hkI} devoted to $\ZZ^2$,
a \emph{proper internal pattern} $\sP$  is defined as a finite SAW with the property
that, for any $k \ge 1$, there exists a SAW containing at least $k$ translates
of $\sP$. The pattern theorem states that: for a given proper internal pattern
$\sP$, there exists $a>0$ such that the number of $n$-step SAWs 
from the origin $0$,  
containing fewer than $an$ translates of $\sP$, 
is exponentially smaller than the total $\s_n := \s_n(0)$. 

The lattice $\ZZ^2$ is bipartite, in that its vertices can be coloured black or white 
in such a way that every edge links a black vertex and a white vertex.
The pattern theorem may be used to prove for this bipartite graph
that
\begin{equation*}
\lim_{n\rightarrow\infty}\frac{\sigma_{n+2}}{\sigma_n}=\mu^2.
\end{equation*} 
The proof is based on a surgery of SAWs that preserves the parity of their lengths.
The following stronger statement has been open since Kesten's paper \cite{hkI},
see the discussion at \cite[p.\ 244]{ms}.

\begin{question}
Is it the case for SAWs on $\ZZ^2$ that $\s_{n+1}/\s_n \to \mu$?
\end{question}

Hammersley's Theorem \ref{jmh} establishes the existence of the 
connective constant for any infinite 
\emph{quasi-transitive} graph. It is easy to construct examples of 
(non-quasi-transitive) graphs for which the limit 
defining $\mu$ does not exist, and it is natural to enquire of the situation for a random graph. For concreteness,
we consider here the infinite cluster $I$ of bond percolation on $\ZZ^2$ with edge-density $p>\frac12$
(see \cite{G99}). 

\begin{question}
Let $\si_n(v)$ be the number of $n$-step SAWs on $I$ starting
at the vertex $v$. Does the limit $\mu(v):= \lim_{n\to\oo} \si_n(v)^{1/n}$
exist a.s., and satisfy $\mu(v)=\mu(w)$ a.s.\ on the event $\{v,w \in I\}$?
\end{question}

Discussions of issues around this question, including of when $\mu(v)=p\mu(\ZZ^2)$ a.s.\ 
on the event $\{v \in I\}$, may be found in papers of Lacoin \cite{Lac1, Lac2}.
The SAW problem on (deterministic) weighted graphs is considered in \cite{GL-wtsaw}
(see Section \ref{ssec:wCg}). 

\subsection{Critical exponents for SAWs}\label{ssec:ce}
\lq Critical exponents' play a significant role in the theory of phase transitions.
Such exponents have natural definitions for SAWs on a given graph, as summarised next. 
The reader is referred to \cite{bdgs, ms} and the references therein
for general accounts of critical exponents for SAWs.
The three exponents that have received most attention in the study of SAWs are 
as follows.

We consider only the case of SAWs in Euclidean spaces, thus excluding,
for example, the hyperbolic space of \cite{MadW}.
Suppose for concreteness that there 
exists a periodic, locally finite embedding of $G$ into $\RR^d$ with $d\geq 2$, and no
such embedding into $\RR^{d-1}$. The case of 
general $G$ has not been studied extensively, and most attention has been paid to the
hypercubic lattice $\ZZ^d$. 

\subsubsection*{The critical exponent $\g$.}
It is believed (when $d \ne 4$) that the generic behaviour of
$\s_n(v)$ is given by:
\begin{equation}\label{ex1}
\sigma_n(v)\sim A_vn^{\g-1}\mu^n,\qquad \text{as}\ n\to\infty,\text{ for }  v\in V,
\end{equation}
for constants $A_v>0$ and  $\g\in\RR$.
The value of the `critical exponent' $\g$ is believed to depend on $d$ only, 
and not further on the choice of
graph $G$. Furthermore, it is believed (and largely proved, see the account in
\cite{ms}) that
$\g = 1$ when $d \ge 4$.
In the borderline case $d=4$, \eqref{ex1} should hold with $\g=1$ and subject to  
the correction factor $(\log n)^{1/4}$. (See the related work \cite{BBS15} 
on weakly self-avoiding walk.)

\subsubsection*{The critical exponent $\eta$.}
Let $v,w\in V$, and
\begin{equation*}
Z_{v,w}(x)=\sum_{n=0}^{\infty}\sigma_n(v,w)x^n,\qquad x>0,
\end{equation*}
where $\sigma_n(v,w)$ is the number of $n$-step SAWs with endpoints $v,w$. It is known under certain 
circumstances that the generating functions $Z_{v,w}$ have radius of convergence $\mu^{-1}$ (see \cite[Cor.\ 3.2.6]{ms}),
and it is believed that there exists an exponent $\eta$ and constants $A_v'>0$ such that
\begin{equation}\label{eq:defeta}
Z_{v,w}(\mu^{-1})\sim A_v'd_G(v,w)^{-(d-2+\eta)},\qquad \text{as } d_G(v,w)\to\infty,
\end{equation}
where $d_G(v,w)$ is the graph-distance between $v$ and $w$.
Furthermore, $\eta$ satisfies $\eta=0$ when $d \ge 4$.  

\subsubsection*{The critical exponent $\nu$.}

Let $\Sigma_n(v)$ be the set of $n$-step SAWs from $v$, and let $\pi_n$ be chosen
at random from $\Si_n(v)$ according to the uniform probability measure.
Let $\|\pi\|$ be the 
graph-distance between the endpoints of a SAW $\pi$. 
It is believed (when $d\neq 4$) that there exists an
exponent $\nu$ (the so-called \emph{Flory exponent}) and constants $A_v''>0$, such that
\begin{equation}\label{eq:defnu}
\EE(\|\pi_n\|^2) \sim A_v''n^{2\nu}, \qquad v\in V.
\end{equation}
As above, this should hold for $d=4$ subject to the inclusion of the 
correction factor $(\log n)^{1/4}$. It is believed that
$\nu=\frac12$ when $d \ge 4$.

The exponent $\nu$ is an indicator of the geometry of an $n$-step SAW $\pi$ 
chosen with the uniform measure.
In the \emph{diffusive} case, we have $\nu=\frac12$, whereas in the \emph{ballistic} case (with
$\|\pi_n\|$ typically of order $n$), we have $\nu=1$.
We return to this exponent in Section \ref{sec:nu}.
 
\smallskip

The three exponents $\g$, $\eta$, $\nu$
are believed to be related through the so-called \emph{Fisher relation}
$\gamma=\nu(2-\eta)$. The definitions \eqref{ex1}, \eqref{eq:defeta}, \eqref{eq:defnu}
may be weakened to logarithmic asymptotics, in which case we say they hold \emph{logarithmically}.

\section{Bounds for connective constants}
\label{sec:lb}

We discuss upper and lower bounds for connective constants in this section, 
beginning with some 
algebraic background.

\subsection{Transitivity of graphs}\label{ssec:trans}
The automorphism group of the graph $G=(V,E)$ is denoted $\Aut(G)$,
and the identity automorphism is written $\id$. 
The expression $\sA\le\sB$ means that $\sA$ is a subgroup of $\sB$, and 
$\sA\normal\sB$ means that $\sA$ is a normal subgroup.

A subgroup $\Ga\le\Aut(G)$ is said to \emph{act transitively} on 
$G$ if, for $v,w\in V$,
there exists $\g\in\Ga$ with $\g v=w$. It is said to \emph{act quasi-transitively} 
if there exists a finite set $W$
of vertices such that, for $v\in V$, there exist 
$w\in W$ and $\g\in \Ga$ with
$\g v=w$. The graph is called \emph{transitive} 
(respectively, \emph{quasi-transitive}) 
if $\Aut(G)$ acts transitively
(respectively, quasi-transitively) on $G$.

An automorphism $\g$ is said to \emph{fix} a vertex $v$ if $\g v = v$.
The \emph{stabilizer} of $v \in V$ is the subgroup
$$
\Stab_v:=\{\g\in\Aut(G): \g v=v\}.
$$
The subgroup $\Ga$ is said to \emph{act freely} on $G$ (or on the vertex-set $V$) if
$\Ga \cap \Stab_v =\{\id\}$ for $v \in V$.

Let $\sG$ (\resp, $\sQ$) be the set of infinite, simple, locally finite, transitive
(\resp, quasi-transitive), rooted graphs, and
let $\sG_\De$ (\resp, $\sQ_\De$) be the subset comprising $\De$-regular graphs. We write
$\id=\id_G$ for the root of the graph $G$.

\subsection{Bounds for $\mu$ in terms of degree}
Let $G$ be an infinite, connected, $\De$-regular graph.
How large or small can $\mu(G)$ be?  It is trivial by counting non-backtracking walks that
$\s_n(v) \le \De(\De-1)^{n-1}$, whence $\mu(G) \le \De -1$
with equality if $G$ is the $\De$-regular tree.
It is not difficult to prove the strict inequality
$\mu(G) < \De-1$ when $G$ is quasi-transitive and contains a cycle
(see \cite[Thm 4.2]{GrL1}).
Lower bounds are harder to obtain.

A multigraph is called \emph{loopless} if each edge has distinct endvertices.

\begin{theorem} {\cite[Thm 4.1]{GrL1}} \label{thm:lower}
Let $G$ be an infinite, connected, $\De$-regular, transitive, loopless multigraph with $\De\ge 2$.
Then $\mu(G)\geq\sqrt{\De-1}$ if either
\begin{letlist}
\item $G$ is simple, or
\item $G$ is non-simple and $\De \le 4$.
\end{letlist}
\end{theorem}

Note that, for the (non-simple) bridge graph $\BB_\De$ with $\De \ge 2$, we have the 
equality $\mu(\BB_\De)= \sqrt{\De-1}$. 

Here is an outline of the proof of Theorem \ref{thm:lower}.
A SAW is called \emph{forward-extendable} if it is the initial
segment of some infinite SAW. 
Let $\siF_n(v)$ be the number of forward-extendable $n$-step SAWs
starting at $v$.
Theorem \ref{thm:lower} is proved by showing as follows that 
\begin{equation}\label{g1}
\siF_{2n}(v) \ge (\De-1)^n.
\end{equation}
Let $\pi$ be a (finite) SAW from $v$, with final endpoint $w$. For a vertex $x \in \pi$ 
satisfying $x \ne w$, and an edge 
$e\notin\pi$ incident to $x$, the pair $(x,e)$
is called $\pi$-\emph{extendable}
if there exists an infinite SAW starting at $v$ whose initial segment
traverses $\pi$ until $x$, and then traverses $e$. 

First, it is proved subject to a 
certain condition $\Pi$ that, for any $2n$-step forward-extendable SAW  $\pi$,
there are at least $n(\Delta-2)$ $\pi$-extendable pairs. 
Inequality \eqref{g1} may be deduced from this statement.

The second part of the proof is to show that graphs
satisfying either (a) or (b) of the theorem satisfy 
condition $\Pi$. It is fairly simple to
show that (b) suffices, and it may well be reasonable to
extend the conclusion to include values of $\De \ge 5$.

\begin{question}\label{qn:q4}
Is it the case that $\mu(G) \ge \sqrt{\De-1}$
in the non-simple case of Theorem \ref{thm:lower}(b) with $\De \ge 5$?
\end{question}

The growth rate $\muF$ of the number of forward-extendable SAWs has been studied 
further by Grimmett, Holroyd, and Peres \cite{GHLP}.
They show that $\muF = \mu$ for any infinite, connected, quasi-transitive
graph, with further results involving the numbers of
backward-extendable and doubly-extendable SAWs. 

\begin{question}
Let $\De\ge 3$.
What is the sharp lower bound $\mumin(\De):=\inf\{\mu(G): G \in \sG_\De\}$?
How does $\mumin(\De)$ behave as $\De\to\oo$?
\end{question}

This question is considered in Section \ref{ssec:cubic} when $\De=3$, and it is asked
in Question \ref{qn:gmean} whether or not $\mumin(3)=\phi$, the golden mean. The lower bound
$\mu\ge \sqrt{\De-1}$ of Theorem \ref{thm:lower}(a)
may be improved as follows when $G$ is non-amenable.

Let $P$ be the transition matrix of simple random walk on $G=(V,E)$, and let $I$ be the 
identity matrix. The \emph{spectral bottom} of $I-P$ is defined to be the largest  
$\lambda$ with the property that, for all $f\in\ell^2(V)$, 
\begin{equation}\label{eq:spbo}
\bigl\langle f, (I-P)f \bigr\rangle\geq \lambda\langle f,f\rangle.
\end{equation}
It may be seen that $\lambda(G)=1-\rho(G)$ where $\rho(G)$ is
the spectral radius of $P$ (see \cite[Sect.\ 6]{LyP}, and \cite{Wo} for
an account of the spectral radius). It is known that $G$ is a non-amenable if and only if $\rho(G)<1$,
which is equivalent to $\lambda(G)>0$. This was proved
by Kesten \cite{K59b, K59a} for Cayley graphs of finitely-presented groups, and 
extended to general transitive graphs by  Dodziuk \cite{Dod}
(see also the references in \cite[Sect.\ 6.10]{LyP}).

\begin{theorem}\cite[Thm 6.2]{GrL7}\label{thm:newin}
Let $G\in \sG_\De$ with $\De \ge 3$, and let $\lambda=\lambda(G)$ be
the above spectral bottom. The connective constant satisfies 
\begin{equation}\label{eq:newin}
\mu(G)\geq (\De-1)^{\frac12(1+c\lambda)},
\end{equation}
where $c=\De(\De-1)/(\De-2)^2$.
\end{theorem}

\subsection{Upper bounds for $\mu$ in terms of degree and girth}

The \emph{girth} of a simple graph is the length of its shortest
cycle. Let $\sG_{\De,g}$ be the subset of $\sG_\De$ containing graphs with girth $g$.

\begin{theorem}\cite[Thm 7.4]{GrL6} \label{g5}
For $G \in \sG_{\De,g}$ where $\De,g \ge 3$, we have that
$\mu(G) \le y$ where $\zeta := 1/y$ is the smallest positive real  root of the equation 
\begin{equation}\label{eq:theeq}
(\De-2)\frac {M_1(\zeta)}{1+M_1(\zeta)} + \frac{M_2(\zeta)}{1+M_2(\zeta)}=1,
\end{equation}
with
\begin{equation} \label{eq:theeq2}
M_1(\zeta)=\zeta, \qquad M_2(\zeta)=2(\zeta+\zeta^2+\dots + \zeta^{g-1}).
\end{equation}
The upper bound $y$ is sharp, and is achieved by the free product graph 
$F := K_2 * K_2 * \dots * K_2 * \ZZ_g$, with $\De-2$ copies
of the complete graph $K_2$ on two vertices and one copy
of the cycle $\ZZ_g$ of length $g$.
\end{theorem}

The proof follows quickly by earlier results of Woess \cite{Wo}, and Gilch and M\"uller \cite{Gilch}.
By \cite[Thm 11.6]{Wo}, every $G \in \sG_{\De,g}$ is covered by $F$, and by 
\cite[Thm 3.3]{Gilch},
$F$ has connective constant $1/\zeta$.

\section{Cubic graphs and the golden mean}\label{sec:fisher}

A graph is called \emph{cubic} if it is regular with degree $\De=3$. 
Cubic graphs have the property that every edge-self-avoiding cycle is also vertex-self-avoiding.
We assume throughout this section that $G=(V,E)\in\sQ_3$,
and we write $\phi:= \frac12(1+\sqrt 5)$ for the golden mean.

\subsection{The Fisher transformation}\label{ssec:fisher}

\begin{figure}[htb]
 \centering
    \includegraphics[width=0.5\textwidth]{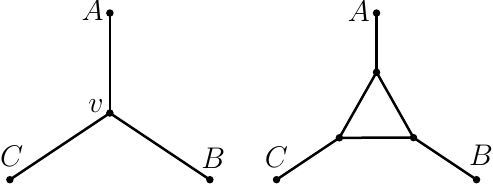}
  \caption{The Fisher transformation of the star.}
  \label{fig:fisher}
\end{figure}

Let $v\in V$, and recall that $v$ has degree $3$ by assumption. 
The so-called \emph{Fisher transformation} acts at $v$ by 
replacing it by a triangle, as illustrated in 
Figure \ref{fig:fisher}. 
The Fisher transformation has been valuable in the study of the relations 
between Ising, dimer, and general vertex models
(see \cite{bdet,fisher,zli,li}), and also
in the calculation of the connective constant of the Archimedean lattice $(3,12^2)$ (see,
for example, \cite{g,GPR,jg}).
The Fisher transformation may be applied at every vertex of a cubic 
graph,
of which the hexagonal and square/octagon lattices are examples.

Let $G$ be, in addition, quasi-transitive. By Theorem \ref{jmh}, $G$ has
a well-defined connective constant $\mu=\mu(G)$ satisfying \eqref{connconst}.
Write $F(G)$ for the graph obtained by applying the Fisher transformation
at every vertex of $G$. 
The automorphism group of $G$ induces an automorphism subgroup of $F(G)$,
so that $F(G)$ is quasi-transitive and has a well-defined connective constant.
It is noted in \cite{GrL2}, and probably elsewhere also, that the connective constants of
$G$ and $F(G)$ have a simple relationship. This conclusion, and its iteration, 
are given in the next theorem

\begin{theorem} \cite[Thm 3.1]{GrL2} \label{thm:main2}
Let $G\in\sQ_3$, and consider the
sequence $(G_k: k=0,1,2,\dots)$ given by $G_0=G$ and $G_{k+1} = \fish(G_k)$.
\begin{letlist}
\item The connective constants $\mu_k := \mu(G_k)$ satisfy
$\mu_k^{-1} = g(\mu_{k+1}^{-1})$ where
$g(x)= x^2 + x^3$.
\item 
The sequence $\mu_k$ converges monotonely to the golden mean $\phi$,
and 
\begin{equation*}
- \left(\frac 47\right)^k \le \mu_k^{-1} - \phi^{-1} 
\le \left(\frac2{7-\sqrt 5}\right)^k, \qquad k \ge 1.
\end{equation*}
\end{letlist}
\end{theorem}

The idea underlying part (a) is that, at each vertex $v$ visited by a SAW $\pi$ on $G_k$, one 
may replace that vertex by either of the two paths around the `Fisher triangle' of $G_{k+1}$
at $v$. Some book-keeping is necessary with this argument, and this is best done 
via the generating functions \eqref{eq:gf}. 

A similar argument may be applied in the context of a `semi-cubic' graph.

\begin{theorem} \cite[Thm 3.3]{GrL2} \label{sf}
Let $G$ be an infinite, connected, bipartite graph with vertex-sets coloured black
and white, and suppose the coloured graph is quasi-transitive, and every black vertex has degree $3$. Let
$\wt{G}$ be the graph obtained by applying the Fisher transformation at each black vertex. 
The connective 
constants $\mu$ and $\wt\mu$ of $G$ and $\wt{G}$, respectively, 
satisfy $\mu^{-2}=h(\wt{\mu}^{-1})$, where $h(x)=x^3+x^4$.
\end{theorem}

\begin{figure}[tb]
 \centering
    \includegraphics[width=0.5\textwidth]{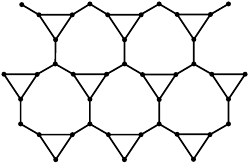}
  \caption{The lattice $\wt\HH$ is derived from the hexagonal lattice $\HH$ by
applying the Fisher transformation at alternate vertices. Its connective constant
$\wt\mu$ is a root of the equation
$x^{-3}+x^{-4} = 1/(2+\sqrt 2)$.}
  \label{fig:newlatt}
\end{figure}

\begin{example}
Take $G=\HH$, the hexagonal lattice with connective constant 
$\mu=\sqrt{2+\sqrt{2}}\approx 1.84776$, see \cite{ds}. The ensuing lattice $\wt{\HH}$ is illustrated
in Figure \ref{fig:newlatt}, and its connective constant $\wt{\mu}$ satisfies $\mu^{-2}=h(\wt{\mu}^{-1})$,
which may be solved to obtain $\wt{\mu}\approx 1.75056$.
\end{example}

We return briefly to the  critical exponents of Section \ref{ssec:ce}.
In \cite[Sect.\ 3]{GrL2}, reasonable definitions of the three exponents 
$\g$, $\eta$, $\nu$ are presented,
none of which depend on the existence of embeddings into $\RR^d$. 
Furthermore, it is proved that the values of the
exponents are unchanged under the Fisher transformation.

\subsection{Bounds for connective constants of cubic graphs}\label{ssec:cubic}

Amongst cubic graphs,  the $3$-regular tree $T_3$ has largest connective constant $\mu(T_3)=2$.
It is an open problem to determine the sharp lower bound on $\mu(G)$ for $G \in \sG_3$.
Recall the ladder graph $\LL$ of Figure \ref{fig:ladder-hex}, with $\mu(\LL)=\phi$.

\begin{question}\cite[Qn 1.1]{GrL6}\label{qn:gmean}
Is it the case that $\mu(G) \ge \phi$ for $G \in \sG_3$?
\end{question}

Even in the case of graphs with small girth, the best general lower bounds 
known so far are as follows.

\begin{theorem}\cite[Thms 7.1, 7.2]{GrL6}\label{g3}
\begin{letlist}
\item For $G\in\sG_{3,3}$, we have that
\begin{equation}\label{uc-2}
\mu(G) \ge x,
\end{equation}
where $x\in(1,2)$ satisfies
\begin{align}
\frac{1}{x^2}+\frac{1}{x^3}&=\frac{1}{\sqrt{2}}.\label{up1}
\end{align}
\item
For $G\in \sG_{3,4}$, we have that
\begin{equation}\label{uc-1}
\mu(G) \ge 12^{1/6}.
\end{equation}
\end{letlist}
\end{theorem}

The sharp upper bounds for $\sG_{3,3}$ and $\sG_{3,4}$
are those of Theorem \ref{g5}, and they are attained
respectively by the Fisher graph of the $3$-regular tree, and of the
degree-$4$ tree (in which each vertex is replaced by a $4$-cycle).

Question \ref{qn:gmean} is known to have a positive answer for various classes of graph, including
so-called TLF-planar graphs (see \cite{GrL6,DR09}).
The word \emph{plane} means a simply connected Riemann surface 
without boundaries. An \emph{embedding} of a graph $G=(V,E)$ in a plane $\sP$ is a function
$\eta:V\cup E \to \sP$ such that $\eta$ restricted to $V$
is an injection and, for $e=\langle u,v\rangle\in E$, $\eta(e)$ is a $C^1$ image of $[0,1]$. 
An embedding is ($\sP$-)\emph{planar} if the images of
distinct edges are disjoint except possibly at their endpoints, 
and a graph is ($\sP$-)\emph{planar} if
it possesses a ($\sP$-)planar embedding.
An embedding is  \emph{topologically locally finite} (\emph{TLF}) if the images of
the vertices have no accumulation point, and a
connected graph is called \emph{TLF-planar} if it possesses a planar TLF embedding. 
Let $\TLF_\De$ denote the class
of transitive,  TLF-planar graphs with vertex-degree $\De$. 

\begin{theorem}\cite{GrL6}\label{tm93}
Let $G\in\TLF_3$ be infinite. Then $\mu(G)\geq \phi$.
\end{theorem}

Two techniques are used repeatedly in the proof. The first is to construct an injection
from the set of  SAWs on the ladder graph $\LL$ that move either rightwards or vertically,
into the set of SAWs on a graph $G'$ derived from $G$. A large subclass of $\sT_3$ may be
treated using such a construction.  The remaining graphs in $\sT_3$ require
detailed analyses using a variety of transformations
of graphs including the Fisher transformation of Section \ref{ssec:fisher}.

A second class of graphs for which $\mu\ge \phi$ is given as follows. 
The definition of a transitive \ghf\ is deferred 
to Definition \ref{def:height}.

\begin{theorem}\cite[Thm 3.1]{GrL6}
We have that $\mu(G) \ge \phi$ for any cubic graph $G \in \sG_3$ that
possesses a transitive \ghf.
\end{theorem}

This theorem covers all Cayley graphs of finitely presented groups with strictly positive 
first Betti numbers (see Section \ref{ssec:cayley} and \cite[Example 3]{GrL6}). 
Cayley graphs are introduced in
Section \ref{ssec:cayley}.

\section{Strict inequalities for connective constants}\label{sec:si}

\subsection{Outline of results}
Consider a probabilistic model on a graph $G$, such as the percolation or random-cluster model 
(see \cite{G-rcm}). There is a parameter (perhaps `density' $p$ or `temperature' $T$)
and a `critical point' (usually written $\pc$ or $\Tc$). The numerical value
of the critical point depends  on the choice of graph $G$.
It is often important to understand whether a systematic change in the graph
causes a \emph{strict} change in the value of the critical point.
A general approach to this issue was presented by Aizenman and Grimmett \cite{AG} 
and developed further in \cite{BBR, BGK,G94} and \cite[Chap.\ 3]{G99}.
The purpose of this section is to review work of \cite{GrL3} directed
at the corresponding question for self-avoiding walks.

Let $G$ be a subgraph of $G'$, and suppose each graph is quasi-transitive.
It is trivial that $\mu(G) \le \mu(G')$. Under what conditions does 
the strict inequality $\mu(G) < \mu(G')$ hold? Two sufficient
conditions for the strict inequality 
are reviewed here.
This is followed in Section \ref{sec:cayley} with a summary 
of  consequences for Cayley graphs.

The results of this section apply to \emph{transitive} graphs. 
Difficulties arise under the weaker assumption of quasi-transitivity.

\subsection{Quotient graphs}\label{sec:quot}
Let $G=(V,E)\in\sG$. 
Let $\Ga\le\Aut(G)$ act transitively, and let $\sA\normal\Ga$
(we shall discuss the non-normal case later). There are several ways
of constructing a quotient graph $G/\sA$, the strongest of which (for our purposes) is
given next.
The set of neighbours of a vertex $v\in V$ is denoted by $\pd v$.

We denote by $\vG =(\olV,\vE)$ the \emph{directed} quotient graph
$G/\sA$ constructed as follows. Let $\approx$ be the equivalence relation on $V$ given by
$v_1 \approx v_2$ if and only if there exists
$\a\in\sA$ with $\a v_1=v_2$. The vertex-set $\olV$
comprises the equivalence classes of $(V,\approx)$, that is, the orbits $\ol v := \sA v$
as $v$ ranges over $V$. For $v,w \in V$, we place $|\pd v \cap \ol w|$ directed edges
from $\ol v$ to $\ol w$ (if $\ol v = \ol w$,
these edges are directed loops).

\begin{example}
Let $G$ be the square lattice $\ZZ^2$ and let $m \ge 1$.
Let $\Ga$ be the set of translations of $\ZZ^2$, and let $\sA$ be the 
normal subgroup of $\Ga$ generated by the map that sends $(i,j)$ to
$(i+m,j)$. The quotient graph $G/\sA$ is the square lattice `wrapped
around a cylinder', with each edge replaced by two oppositely directed edges.
\end{example}

Since $\vG$ is obtained from $G$
by a process of identification of vertices and edges, it is natural to ask whether 
the strict inequality $\mu(\vG) < \mu(G)$
is valid. Sufficient conditions for this strict inequality are presented next. 

Let $L=L(G,\sA)$ be the length of the shortest SAW of $G$
with (distinct) endpoints in the same orbit.
Thus, for example, $L=1$ if $\vG$ possesses a directed loop.
A group is called \emph{trivial} if it comprises the identity only.

\begin{theorem} \cite[Thm 3.8]{GrL3} \label{si}
Let $\Ga$ act transitively on $G$, and let
$\sA$ be a non-trivial, normal subgroup of $\Ga$. 
The connective constant $\vmu=\mu(\vG)$ satisfies $\vmu < \mu(G)$ if: either
\begin{letlist}
\item $L \ne 2$, or
\item $L=2$ and either of the following holds:
\begin{romlist}
\item $G$ contains some $2$-step
 SAW $v\ (=w_0),w_1,w_2\ (=v')$ satisfying $\ol v = \ol v'$ and $|\pd v \cap \ol w_1| \ge 2$, 
\item $G$ contains some SAW $v\,(=w_0),w_1,w_2,\dots,w_l\,(=v')$ satisfying 
$\ol v = \ol v'$, $\ol w_i \ne \ol w_j$ for $0 \le i < j < l$, and furthermore $v' = \a v$ for
some $\a \in \sA$ which fixes no $w_i$.
\end{romlist}
\end{letlist}
\end{theorem} 

\begin{remark}\label{qn:strict}
In the situation of Theorem \ref{si}, can one \emph{calculate} an explicit $R=R(G,\sA)<1$
such that $\mu(\vG)/\mu(G) <R$? The answer is (in principle)
positive under a certain condition, namely
that the so-called \lq bridge constant' of $G$ equals its connective constant.
Bridges are discussed in Section \ref{sec:bridge}, and it is shown in
Theorem \ref{thm:bridge} that the above holds when $G$ possesses a so-called \lq unimodular \ghf'
(see Definition \ref{def:height}). See also \cite[Thm 3.11]{GrL3} and \cite[Remark 4.5]{GrL4}.
\end{remark}

We call $\sA$ \emph{symmetric} if 
\begin{equation*}
|\partial v\cap\ol{w}|=|\partial w\cap\ol{v}|,\qquad v,w\in V.
\end{equation*}
Consider the special case $L=2$ of Theorem \ref{si}.
Condition (i) of Theorem \ref{si}(b)
holds if $\sA$ is symmetric, since $|\pd w \cap \ol v| \ge 2$.
Symmetry of $\sA$ is implied by unimodularity, for a definition
of which we refer the 
reader to \cite[Sect.\ 3.5]{GrL3} or \cite[Sect.\ 8.2]{LyP}.

\begin{example}
Conditions (i)--(ii) of Theorem \ref{si}(b) are necessary in the case $L=2$,
in the sense illustrated by the following example. Let $G$ be the 
infinite $3$-regular tree with a distinguished end $\om$.
Let $\Ga$ be the set of automorphisms
that preserve $\omega$, and let $\sA$ be the normal subgroup
generated by the interchanges of the  two children of any given vertex $v$
(and the associated relabelling of their descendants). 
The graph $\vG$ is isomorphic to that obtained from $\ZZ$ by
replacing each edge by two directed edges in one direction and one in the reverse direction.
It is easily seen that $L=2$, but that
neither (i) nor (ii) holds.
Indeed, $\mu(\vG) = \mu(G) = 2$.
\end{example}

The proof of Theorem \ref{si} follows partly the general approach of Kesten
in his pattern theorem, see \cite{hkI} and \cite[Sect.\ 7.2]{ms}.
Any $n$-step
SAW $\vpi$ in the directed graph $\vG$  lifts to a SAW $\pi$ in the larger
graph $G$. The
idea is to show there exists $a>0$ such that `most' such $\vpi$ contain at least $an$ sub-SAWs 
for which the corresponding sub-walks of $\pi$ may be replaced by SAWs
on $G$.
Different subsets of these sub-SAWs of $\vG$ give rise to different SAWs 
on $G$. The number of such subsets grows exponentially in $n$, and this
introduces an exponential `entropic' factor in the counts of SAWs.

Unlike Kesten's proof and its later elaborations, these results
apply in the general setting of transitive graphs, and they
utilize algebraic and combinatorial techniques.

We discuss next the assumption of normality of $\sA$ in
Theorem \ref{si}. The (undirected) simple quotient graph $\olG=(\olV,\olE)$ may be defined as follows
even if $\sA$ is not a normal subgroup of $\Ga$. As before, the vertex-set $\olV$ is the set
of orbits of $V$ under $\sA$. Two distinct orbits $\sA v$, $\sA w$ are declared adjacent 
in $\olG$ if there exist $v' \in \sA v$ and $w' \in \sA w$ with $\langle v',w'\rangle \in E$.
We write $\olG = G_\sA$ to emphasize the role of $\sA$.

The relationship between the \emph{site percolation} critical points of $G$ and 
$G_\sA$ is the topic of
a conjecture of Benjamini and Schramm \cite{BenS96}, which
appears to make the additional assumption that $\sA$
acts freely on $V$.  The last assumption is stronger
than the assumption of unimodularity.

We ask for an example in which the non-normal case is essentially different from the
normal case.

\begin{question}\label{qn:q5}
Let $\Ga$ be a subgroup of $\Aut(G)$ acting transitively on $G$. 
Can there exist a non-normal subgroup $\sA$
of $\Ga$ such that: {\rm(i)} the quotient graph $G_\sA$ is transitive,
and {\rm(ii)} there exists no normal subgroup $\sN$ of some
transitively acting $\Ga'$ such that
$G_\sA$ is isomorphic to $G_\sN$? 
Might it be relevant to assume that $\sA$ acts freely on $V$?
\end{question}

We return to connective constants with the following question, inspired in
part by \cite{BenS96}.

\begin{question}\label{qn:qn7'}
Is it the case that $\mu(G_\sA) < \mu(G)$ under the assumption that 
$\sA$ is a non-trivial (not necessarily normal) subgroup of $\Ga$ 
acting freely on $V$, such that $G_\sA$ is transitive?
\end{question}

The proof of Theorem \ref{si}  may be adapted
to give an affirmative answer to Question \ref{qn:qn7'} subject 
to a certain extra condition on $\sA$,
see \cite[Thm 3.12]{GrL3}. Namely, it suffices that there exists $l \in \NN$
such that $G_\sA$ possesses a cycle of length $l$ but $G$ has no cycle of this length.

\subsection{Quasi-transitive augmentations} \label{sec:qta}
We consider next the systematic addition of new edges, and the effect thereof
on the connective constant.
Let $G=(V,E) \in \sG$. From $G$, we derive a second graph
$\ollG=(V,\ollE)$ 
by adding further edges to $E$, possibly in parallel to existing edges. 
\emph{We assume that $E$ is a proper subset of $\ollE$}.

\begin{theorem} \cite[Thm 3.2]{GrL3}\label{qta}
Let $\Ga\le\Aut(G)$ act transitively on $G$, and 
let $\sA\le \Ga$ satisfy either or both of the following.
\begin{letlist}
\item $\sA$ is a normal subgroup of $\Ga$ acting quasi-transitively on $G$.
\item The index $[\Ga:A]$ is finite.
\end{letlist}
If $\sA\le \Aut(G')$, then $\mu(G)<\mu(\ollG)$.
\end{theorem}

\begin{example} \label{ex:sqt}
Let $\ZZ^2$ be the square lattice, with $\sA$  the group
of its translations.  The triangular lattice $\TT$ is obtained from $\ZZ^2$ by adding the edge 
$e=\langle 0,(1,1)\rangle$
together with its images under $\sA$, where $0$ denotes the origin. 
Since $\sA$ is a normal subgroup of itself, it follows that $\mu(\ZZ^2) < \mu(\TT)$. This example
may be extended to augmentations by other periodic families of new edges,
as explained in \cite[Example 3.4]{GrL3}.
\end{example}

\begin{remark}\label{qn:strict2}
In the situation of Theorem \ref{qta}, 
can one calculate an $R>1$ such that $\mu(\ollG)/\mu(G) > R$?
As in Remark \ref{qn:strict}, the answer is positive when $G'$ has a unimodular
\ghf.
\end{remark}

A slightly more general form of Theorem \ref{qta} is presented in \cite{GrL3}.
Can one dispense with the assumption of normality in Theorem \ref{qta}(a)?

\begin{question}\label{qn:nonnormal}
Let $\Ga$ act transitively on $G$, and let $\sA$ be a subgroup of $\Ga$
that acts quasi-transitively on $G$. 
If $\sA \le \Aut(\ollG)$, is it necessarily the case that
$\mu(G) < \mu(\ollG)$?
\end{question}

A positive answer would be implied by an affirmative answer to the following question.

\begin{question}\label{qn:nonnormal2}
Let $G\in\sG$, and let $\sA\le\Aut(G)$
act quasi-transitively on $G$.  When does there exist 
a subgroup $\Ga$ of $\Aut(G)$ acting transitively on $G$ such that $\sA \le \Ga$ and
$\Ga\cap\Stab_v \le \sA$ for $v \in V$?
\end{question}

See \cite[Prop.\ 3.6]{GrL3} and the further discussion therein.

\section{Connective constants of Cayley graphs}\label{sec:cayley}

\subsection{Cayley graphs}\label{ssec:cayley}
Let $\Ga$ be an infinite group with identity element $\id$ and finite generator-set $S$, 
where $S$ satisfies $S=S^{-1}$ and $\id \notin S$. Thus, $\Ga$ has a presentation as
$\Ga = \langle S \mid R\rangle$ where $R$ is a set of relators.
The group $\Ga$ is called \emph{finitely generated} since $|S|<\oo$, and
\emph{finitely presented} if, in addition, $|R|<\oo$

The Cayley graph $G=G(\Ga,S)$ is defined as follows.
The vertex-set $V$ of $G$ is the set of elements of $\Ga$. 
Distinct elements $g,h \in V$
are connected by an edge if and only if there exists $s \in S$ such that $h=gs$.
It is easily seen that $G$ is connected, and $\Ga$ acts  transitively by
left-multiplication. It is standard that $\Ga$ acts freely, and hence
$G$ is unimodular and therefore symmetric. Accounts of
Cayley graphs may be found in \cite{bab95}  and \cite[Sect.\ 3.4]{LyP}. 

Reference is occasionally made here to the \emph{first Betti number} of $\Ga$.
This is the power of $\ZZ$, denoted $B(\Ga)$, in the abelianization $\Ga/[\Ga,\Ga]$. 
(See  \cite[Remark 4.2]{GrL5}.)

\subsection{Strict inequalities for Cayley graphs}
Theorems \ref{si} and \ref{qta} have the following implications for Cayley graphs.
Let $s_1s_2\cdots s_l = \id$ be a relation. 
This relation corresponds to 
the closed walk $(\id,s_1, s_1s_2,\dots, \break s_1 s_2\cdots s_l=\id)$
of $G$ passing through the identity $\id$. 
Consider now the effect of adding a further relator.
Let $t_1,t_2,\dots,t_l \in S$ be such that 
$\rho := t_1t_2\cdots t_l$ satisfies $\rho \ne \id$, and write
$\Ga_\rho = \langle S \mid R\cup\{\rho\}\rangle$.
Then $\Ga_\rho$ is isomorphic to the quotient group $\Ga/\sN$ where
$\sN$ is the normal subgroup of $\Ga$ generated by $\rho$.

\begin{theorem} \cite[Corollaries 4.1, 4.3]{GrL3}\label{caleydecrease}
Let $G=G(\Ga,S)$ be the Cayley graph of
the infinite, finitely presented group $\Ga = \langle S\mid R\rangle$.

\begin{letlist}
\item
Let $G_\rho= G(\Ga_\rho, S)$ be the Cayley graph obtained by adding to $R$ a further non-trivial relator $\rho$. 
Then $\mu(G_\rho) <\mu(G)$.

\item
Let $w \in\Ga$ satisfy $w \ne \id$, $w \notin S$,
and let $\olG_w$ be the Cayley graph of the group obtained by adding $w$
(and $w^{-1}$) to $S$. Then $\mu(G) < \mu(\olG_w)$.
\end{letlist}
\end{theorem} 

As noted in Remarks \ref{qn:strict} and \ref{qn:strict2}, non-trivial 
bounds may in principle be calculated for the ratios of the two connective constants in case (a) 
(\resp, case (b))
whenever $G$ (\resp, $\olG_w$) has a unimodular \ghf.

\begin{example}\label{solattice}
The \emph{square/octagon} lattice, otherwise known as the Archimedean
lattice $(4,8^2)$, is the 
Cayley graph of the group with generator set 
$S=\{s_1,s_2,s_3\}$ and relator set 
$$
R=\{s_1^2, s_2^2, s_3^2,  s_1s_2s_1s_2,s_1s_3s_2s_3s_1s_3s_2s_3\}.
$$
(See \cite[Fig.\ 3]{GrL3}.)
By adding the further
relator $s_2s_3s_2s_3$, we obtain a graph isomorphic
to the ladder graph of
Figure \ref{fig:ladder-hex}, whose connective constant is the golden mean $\phi$. 

By Theorem \ref{caleydecrease}(a),
the connective constant $\mu$ of the square/octagon lattice is strictly greater than 
$\phi = 1.618\dots$. The best lower bound currently known appears to
be $\mu > 1.804\dots$, see \cite{j04}.
\end{example}

\begin{example}
The square lattice $\ZZ^2$ is the Cayley graph
of the abelian group with $S = \{a,b\}$
and $R=\{aba^{-1}b^{-1}\}$. We add a generator $ab$
(and its inverse), thus adding a diagonal to each square of $\ZZ^2$.
Theorem \ref{caleydecrease}(b) implies the standard inequality $\mu(\ZZ^2)< \mu(\TT)$
of Example \ref{ex:sqt}.
\end{example}

\section{Bridges} \label{sec:bridge}

\subsection{Bridges and \ghf s}
Various surgical constructions are useful in the study 
of self-avoiding walks, of which the most elementary involves concatenations of so-called
\lq bridges'. Bridges were introduced by Hammersley and Welsh \cite{HW62}
in the context of the hypercubic lattice $\ZZ^d$. An \emph{$n$-step bridge} on $\ZZ^d$
is a self-avoiding walk $\pi=(\pi_0,\pi_1, \dots, \pi_n)$ such that
$$
\pi_0(1) < \pi_m(1) \le \pi_n(1), \qq 0 < m \le n,
$$
where $x(1)$ denotes the first coordinate of a vertex $x\in\ZZ^d$. 

The significant property
of bridges is as follows: given two bridges $\pi=(0,x_1,\dots, x_m)$, $\pi'=(0,y_1,\dots,y_n)$ 
starting at $0$, the concatenation $\pi \cup [x_m+\pi']$ is
an $(m+n)$-step bridge from $0$. It follows that the number $b_n$ of $n$-step bridges from $0$ satisfies
\begin{equation}\label{eq:superm}
b_{m+n} \ge b_mb_n,
\end{equation}
whence the \emph{bridge constant} $\beta(\ZZ^d):=\lim_{n\to\oo} b_n^{1/n}$ exists. 
Since $b_n \le \si_n$, it is trivial that $\beta(\ZZ^d)\le \mu(\ZZ^d)$. 
Using a surgery argument, Hammersley and Welsh proved
amongst other things that $\beta=\mu$ for $\ZZ^d$.

In this section, we discuss the bridge constant for transitive graphs, 
therein introducing the \ghf s that will be useful in the discussion of locality
in Section \ref{sec:loc}.

First we define a \ghf, and then we use such a function to define a bridge.

\begin{definition} \cite[Defn 3.1]{GrL4}\label{def:height}
Let $G =(V,E)\in \sQ$ with root labelled $\id$. 
\begin{romlist}
\item
A \emph{\ghf} on $G$ is a pair $(h,\sH)$ such that:
\begin{letlist}
\item $h:V \to\ZZ$, and $h(\id)=0$, 
\item $\sH\le \Aut(G)$ acts quasi-transitively on $G$ 
such that $h$ is \emph{\hdi}, in the sense that
$$
h(\a v) - h(\a u) = h(v) - h(u), \qq \a \in \sH,\ u,v \in V,
$$
\item for  $v\in V$,
there exist $u,w \in \pd v$ such that
$h(u) < h(v) < h(w)$.
\end{letlist}
\item
A \ghf\ $(h,\sH)$ is called \emph{transitive} (\resp, \emph{unimodular}) if the action
of $\sH$ is transitive (\resp, unimodular).
\end{romlist}
\end{definition}

The reader is referred to \cite[Chap.\ 8]{LyP} and \cite[eqn (3.1)]{GrL4} 
for discussions of unimodularity.

\subsection{The bridge constant}
Let $(h,\sH)$ be a \ghf\ of the graph $G\in \sQ$. A \emph{bridge} $\pi=(\pi_0,\pi_1,\dots,\pi_n)$
is a SAW on $G$ satisfying
$$
h(\pi_0)< h(\pi_m) \le h(\pi_n), \qq 0< m \le n.
$$
Let $b_n$ be the number of $n$-step bridges $\pi$ from $\pi_0=\id$.
Using quasi-transitivity, it may be shown (similarly to \eqref{eq:superm}) 
that the limit $\beta=\lim_{n\to\oo} b_n^{1/n}$
exists, and $\beta$ is called the \emph{bridge constant}. 
Note that $\beta$ depends on the choice of \ghf.

The following is proved by an extension of the methods of \cite{HW62}.

\begin{theorem}\cite[Thm 4.3]{GrL4}\label{thm:bridge}
Let $G\in \sQ$ possess a unimodular \ghf\ $(h,\sH)$. The associated bridge
constant  $\beta=\beta(G,h,\sH)$ satisfies $\beta= \mu(G)$.
\end{theorem}

In particular, the value of $\beta$ does not depend on the choice of \emph{unimodular} \ghf.

\subsection{Weighted Cayley graphs}\label{ssec:wCg}
A natural extension of the theory of self-avoiding walks is to \emph{edge-weighted} graphs.
Let $G=(V,E)$ be an infinite graph, and let $\phi:E \to [0,\oo)$.
The \emph{weight} of a SAW $\pi$  traversing
the edges $e_1,e_2,\dots,e_n$ is defined as 
$$
w_\phi(\pi) := \prod_{i=1}^n \phi(e_i).
$$
One may ask about the asymptotic behaviour of the sum of the $w_\phi(\pi)$ over
all SAWs $\pi$ with length $n$ starting at a given vertex. 
The question is more interesting when $G$ is not assumed locally finite, 
since the number $\si_n$ of SAWs from a given vertex may then be infinite.  
Some conditions are needed on the pair $(G,\phi)$, and these are easiest stated when $G$ is a Cayley graph.

Let $\Ga=\langle S\mid R\rangle$ be an infinite, finitely presented group, with a Cayley graph $G$
(we do not assume that $G$ is locally finite). Let $\phi:\Ga \to [0,\oo)$ be such that
\begin{letlist}
\item $\phi(\id)=0$,
\item $\phi$ is \emph{symmetric} in that $\phi(\g)=\phi(\g^{-1})$ for $\g\in\Ga$,
\item $\phi$ is \emph{summable} in that $\sum_{\g\in\Ga}\phi(\g)<\oo$.
\end{letlist}
The aggregate weights of SAWs on $G$ have been studied in \cite{GL-wtsaw}. 
It turns out to be useful to consider a generalized notion of the length $l(\pi)$ of a SAW
$\pi$, and there is an interaction between  $l$ and  
$\phi$.   Subject to certain assumptions (in particular, $\Ga$ is assumed virtually indicable),
it is shown that the bridge and connective constants are equal. 
This yields a continuity theorem for the connective constants of
weighted graphs. 

\section{Locality of connective constants}\label{sec:loc}

\subsection{Locality of critical values}

The locality question for SAWs may be stated as follows: 
for which families of rooted graphs is the value of the connective constant $\mu=\mu(G)$ determined by 
the graph-structure of large bounded 
neighbourhoods of the root of $G$? Similar questions have been asked for other systems including 
the percolation model, see \cite{bnp,mt}.

Let $G \in \sQ$. The \emph{ball} $S_k=S_k(G)$, with radius $k$, is
the rooted subgraph of $G$ induced by the set of its vertices within
distance $k$ of the root $\id$. For $G,G'\in \sQ$, we write
$S_k(G) \simeq S_k(G')$ if there exists a  graph-isomorphism from $S_k(G)$ to
$S_k(G')$ that maps $\id$ to $\id'$.  Let
$$
K(G,G') = \max\bigl\{k: S_k(G) \simeq S_k(G')\bigr\}, \qq G,G' \in \sG,
$$
and $d(G,G') = 2^{-K(G,G')}$.
The corresponding metric space was introduced by Babai \cite{Bab};
see also \cite{BenS01,DL}.

\begin{question}
Under what conditions on $G\in\sQ$ and $\{G_n\} \subseteq \sQ$ is it the case that
$$
\mu(G_n) \to \mu(G) \q\text{if} \q K(G,G_n) \to\oo?
$$
\end{question}

The locality property of 
connective constants turns out to be related in a surprising way to the existence of 
\emph{harmonic} \ghf s (see Theorem \ref{eag}). 

\subsection{Locality theorem}

Let $G \in  \sQ$ support a \ghf\ $(h,\sH)$. There are two associated integers $d$, $r$ 
defined as follows.
Let
\begin{equation}\label{eq:defd}
d=d(h)=\max\bigl\{|h(u)-h(v)|: u,v\in V,\ u \sim v\bigr\}.
\end{equation}
If $\sH$ acts transitively, we set
$r=0$. Assume $\sH$ does not act transitively, and 
let $r=r(h,\sH)$ be the infimum of all $r$ such that the following holds.
Let $o_1,o_2,\dots,o_M$ be representatives of the orbits of $\sH$.
For $i\ne j$, there
exists $v_j \in \sH o_j$ 
such that $h(o_i)<h(v_j)$, and a SAW from $o_i$
to $v_j$, with length $r$ or less, all of whose vertices $x$, other than
its endvertices,  satisfy $h(o_i)<h(x)< h(v_j)$.

For $D\ge 1$ and $R\ge 0$, let $\sQ_{D,R}$ be the subset of
$\sQ$ containing graphs which possess a \ughf\  
$(h, \sH)$ with $d(h)\le D$ and $r(h,\sH) \le R$.

\begin{theorem}[Locality theorem] \cite[Thm 5.1]{GrL4}
\label{thm2}
Let  $G \in \sQ$.
Let $D\ge 1$ and $R \ge 0$, and let 
$G_n \in \sQ_{D,R}$ for $n \ge 1$. If
$K(G,G_n)  \to \oo$ as $n \to\oo$, then
$\mu(G_n) \to \mu(G)$.
\end{theorem}

The rationale of the proof is as follows. Consider for simplicity the case of transitive graphs.
Since $\log\sigma_n$ is a subadditive sequence (see \eqref{eq:subadditive}), 
we have that $\mu \le \sigma_n^{1/n}$ for $n\ge 1$. 
Similarly, by \eqref{eq:superm}, $\log \beta_n$ is
superadditive, so that $\beta \ge b_n^{1/n}$ for $n \ge 1$.
Therefore,
\begin{equation}\label{eq:sandw}
b_n^{1/n} \le \beta \le \mu \le \sigma_n^{1/n}, \qq n\ge 1.
\end{equation}
Now, $b_n$ and $\sigma_n$ depend only on the ball $S_n(G)$.
If $G$ is such that $\beta=\mu$, then their shared value can be approximated,
within any prescribed accuracy, by counts of bridges and SAWs on bounded balls. 
By Theorem \ref{thm:bridge}, this holds if $G$ supports a \ughf.

\subsection{Application to Cayley graphs}\label{ssec:appl}

Let $\Ga=\langle S \mid R\rangle$ be finitely presented with  Cayley graph $G=G(\Ga,S)$.
Let $t \in \Ga$ have infinite order. We present an application
of the Locality Theorem \ref{thm2} to the situation in which a new relator
$t^n$ is added. Let $G_n$ be the Cayley graph of
the group $\Ga_n=\langle S\mid R\cup\{t^n\}\rangle$.

\begin{theorem}\cite[Thm 6.1]{GrL5}\label{torus}
If the first Betti number of $\Ga$ satisfies $B(\Ga) \ge  2$,  then 
 $\mu(G_n) \to \mu(G)$ as $n\to\oo$.
\end{theorem}

\section{Existence of \ghf s}\label{sec:exist}

We saw in Sections \ref{sec:bridge} and \ref{sec:loc} that, 
subject to the existence of certain \ughf s, 
the equality $\beta=\mu$ holds, and a locality result follows. 
In addition, there exists a terminating algorithm for calculating $\mu$ to any 
degree of precision (see \cite{GrL4}).
In this section, we identify certain classes of graphs that possess \ughf s.

For simplicity in the following,
we restrict ourselves to Cayley graphs of finitely generated groups.

\subsection{Elementary amenable groups}

The class $\EG$ of elementary amenable groups was introduced by Day in 1957, \cite{Day57}, as the 
smallest class of groups that contains the set $\EG_0$ of all finite and abelian groups, 
and is closed under the operations 
of taking subgroups, and of forming quotients, extensions, and directed unions.
Day noted that every group in $\EG$ is amenable (see also von Neumann \cite{vN}). Let $\EGF$ be the set of
infinite, finitely generated members of $\EG$. 

\begin{theorem}\cite[Thm 4.1]{GrL7}\label{eag}
Let $\Ga\in\EGF$. 
There exists a normal subgroup $\sH\normal\Ga$
with $[\Ga:\sH]<\oo$ such that any locally finite Cayley graph $G$ of 
$\Ga$ possesses a \ghf\ of the form $(h,\sH)$ which is both unimodular and harmonic.
\end{theorem}

Note that the \ghf\ $(h,\sH)$ of the theorem is \emph{harmonic}. It has a further property,
namely that $\sH\normal \Ga$ has finite index, and $\sH$ acts on $\Ga$ 
by left-multiplication. Such a \ghf\ is called \emph{strong}.

The proof of Theorem \ref{eag} has two stages. Firstly, by a standard algebraic result,
there exist $\sH\normal \Ga$ such that: 
$|\Ga/\sH|<\oo$, and $\sH$ is \emph{indicable} in that
there exists a surjective homomorphism $F: \sH \to \ZZ$.
At the second stage, we consider a random walk on the Cayley graph $G$,
and set $h(\g)=\EE_\g(F(H))$, where $H$ is the first hitting point of $\sH$
viewed as a subset of vertices.
That $h$ is harmonic off $\sH$ is automatic, and on $\sH$ because the action of $\sH$
is unimodular.

The conclusion of Theorem \ref{eag} is in fact valid for the larger class of 
infinite, finitely generated, virtually indicable groups (see \cite[Thm 3.2]{GL-wtsaw}).

\subsection{Graphs with no \ghf}
There exist
transitive graphs possessing no \ghf, and examples 
include the (amenable) Cayley graph of the Grigorchuk
group, and the (non-amenable) Cayley graph of the
Higman group (see \cite[Thms 5.1, 8.1]{GrL7}).
This may be deduced from the next theorem. 

\begin{theorem}\cite[Cor.\ 9.2]{GrL7}\label{thm:noghf} 
Let $\Ga = \langle S \mid R\rangle$ where $|S|<\oo$, 
and let $\Pi$ be the subgroup of permutations of $S$ that 
preserve $\Ga$ up to isomorphism.  Let $G$ be a Cayley graph of
$\Ga$ satisfying $\Stab_1 = \Pi$, where $\Pi$ is viewed as a subgroup of $\Aut(G)$. 
\begin{letlist}
\item If $\Ga$ is a torsion group, then $G$ has no \ghf.
\item Suppose $\Ga$ has no proper, normal subgroup with finite index. If 
$G$ has \ghf\ $(h,\sH)$, then $(h,\Ga)$ is also a \ghf.
\end{letlist}
\end{theorem}

The point is that, when $\Stab_\id=\Pi$,   every automorphism of $G$ is obtained
by a certain composition of an element of $\Ga$ and an element of $\Pi$.
The Grigorchuk group is a torsion group, and part (a) applies.  
The Higman group $\Ga$ satisfies part (b), and is quickly seen to have no \ghf\
of the form $(h,\Ga)$. (A \ghf\ of the form $(h,\Ga)$ is
called a \emph{group height function} in \cite[Sect.\ 4]{GrL5}.)

\section{Speed, and the exponent $\nu$}\label{sec:nu}

For simplicity in this section, we consider only transitive rooted graphs $G$.
Let $\pi_n$ be a random $n$-step SAW from the root of  $G$, 
chosen according to the uniform measure on the set $\Si_n$ of such walks.
What can be said about the graph-distance $\|\pi_n\|$ between the endpoints of $\pi_n$?

We say that SAW on $G$ has \emph{positive speed} if there exist $c,\alpha>0$ such that
$$
\PP(\|\pi_n\| \le cn) \le e^{-\alpha n}, \qq  n \ge 0.
$$ 
While stronger than the natural definition through the requirement of exponential decay to $0$, 
this is a useful definition for the results
of this section. When $G$ is infinite, connected, and quasi-transitive, and SAW on $G$ has  positive speed, 
it is immediate that 
$$
C n^2 \le \EE(\|\pi_n\|^2) \le n^2,
$$
for some $C>0$; thus \eqref{eq:defnu} holds (in a slightly weaker form) with $\nu=1$.

For SAW on $\ZZ^d$ it is known, \cite{CH13}, that SAW does not have positive speed, and that
that $\EE( \|\pi_n\|^2) / n^2 \to 0$  as $n \to\oo$ (cf.\ \eqref{eq:defnu}).
Complementary `delocalization' results have been proved in \cite{CGHM}, for example that,
for $\eps>0$ and large $n$,
$$
\PP(\|\pi_n\|=1) \le n^{-\frac14+\epsilon},
$$
and it is asked there whether
$$
\PP(\|\pi_n\|=x) \le n^{-\frac14+\epsilon}, \qq x \in \ZZ^d.
$$
  
We pose the following question for non-amenable graphs.
  
\begin{question}\cite{CH13}\label{qn:20}
Is it the case that, for any non-amenable Cayley graph $G$ of an infinite, finitely generated group, 
SAW on $G$ has positive speed?
\end{question}

Progress towards this question  may be summarised as follows.
By  bounding the number of SAWs by the number of non-backtracking paths,
Nachmias and Peres \cite[eqn (2.3)]{NP12} have proved that that, for a 
non-amenable, transitive graph $G$ satisfying
\begin{equation}
(\Delta-1)\rho<\mu ,\label{cd12}
\end{equation}
SAW on $G$ has positive speed. 
Here, $\Delta$ is the vertex degree, $\rho$ is the spectral radius of simple random walk on $G$
(see the discussion around \eqref{eq:spbo}), 
and $\mu$ is the connective constant.

It is classical (see, for example, \cite[eqn (1.13)]{G99})  that
$\mu\pc\ge 1$, where $\pc$ is the critical probability of bond percolation on $G$. Inequality \eqref{cd12}
is therefore implied by the stronger inequality
\begin{equation}\label{cd12+}
(\Delta-1)\rho\pc<1.
\end{equation}
These inequalities \eqref{cd12}--\eqref{cd12+} are useful in several settings.
\begin{Alist}
\item \cite{NP12}
There exist $\rho_0<1$ and
$g_0<\oo$ such that, if $G$ is an non-amenable, transitive graph with spectral radius less than $\rho_0$
and girth at least $g_0$, then \eqref{cd12+} holds, and hence SAW on $G$ has positive speed.

\item
\cite{PS00}
Let $S\subset \Ga$ be a finite symmetric generating set of an infinite group $\Ga$, and $\De=|S|$.
Let $S^{(k)}$ be the multiset of cardinality $\De^k$ comprising all elements $g\in \Ga$ with length 
not exceeding $k$ in the word metric given by $S$, each such element  included with multiplicity equal to the
number of such ways of expressing $g$. The set $S^{(k)}$ generates $\Ga$. By \cite[Proof of Thm 2]{PS00},
\begin{equation*}
\rho(G,S^{(k)})\pc(G,S^{(k)})\De^k\to 0 \qq\text{as } k \to\oo,
\end{equation*}
where $(G,S^{(k)})$ denotes the (possibly non-simple) Cayley graph of $\Ga$ with generator-set $S^{(k)}$.
Inequality \eqref{cd12+} follows for sufficiently large $k$

\item
By the argument of \cite[Prop.\ 1]{PS00}, \eqref{cd12+} holds if
\begin{equation*}
\rho<\frac{\De^2+\De-1}{\De^2+(\De-1)^2}.
\end{equation*}
This holds when $\rho<\frac12$, irrespective of $\De$.

\item
Thom \cite{AT13} has shown that,  for any finitely generated, non-amenable group $\Ga$ and any $\epsilon>0$, 
there exists a finite symmetric 
generating set $S$ such that the corresponding Cayley graph $G=G(\Ga,S)$ has  $\rho<\epsilon$. 
By the above, SAW on $G$ has positive speed. 

\item \cite{RS01} Let $i=i(G)$ be the edge-isoperimetric constant of an infinite, transitive graph $G$. 
Since $\pc \le (1+i)^{-1}$
and $\rho^2 \le 1 - (i/\De)^2$ (\cite[Thm 2.1(a)]{BM88}),  inequality \eqref{cd12+} holds whenever
$$
(\De-1)\frac{\sqrt{1-(i/\De)^2}}{1+i} < 1.
$$
It is sufficient that $i/\De > 1/\sqrt 2$.

\item
\cite{DM10}
It is proved that
$$
\rho \le \frac{\sqrt{8\De-16}+3.47}\De,
$$ 
when $G$ is planar. When combined with C above, for example, this implies that SAW has
positive speed on any transitive, planar graph with sufficiently large $\De$.
A related inequality for hyperbolic tesselations is found  in \cite[Thm 7.4]{DM10}.
\end{Alist}

We turn next to graphs embedded in the hyperbolic plane.
It was proved by Madras and Wu \cite{MW05} that SAWs on most regular tilings of the hyperbolic plane
have positive speed. Note that the graphs treated in \cite{MW05} are both (vertex-)transitive and edge-transitive
(unlike the graphs in F above).
It was proved by Benjamini \cite{IB16} that SAWs on the seven regular planar triangulations of 
the hyperbolic plane
have mean displacement bounded beneath by a linear function.

We turn finally to a discussion of the number of ends of a transitive graph.
The number of \emph{ends} of an infinite, connected graph 
$G=(V,E)$ is the supremum over all finite subsets $W\subset V$ 
of the number of infinite components that remain after deletion of $W$. 
An infinite, finitely generated group $\Ga$ is said to have $k$ ends if some locally finite Cayley graph 
(and hence all such Cayley graphs) has $k$ ends. 
Recall from \cite[Prop.\ 6.2]{BM91} that a transitive graph $G$ has $k\in\{1,2,\oo\}$
and, moreover, if $k=2$ (\resp, $k=\oo$) then $G$ is amenable (\resp, non-amenable).

Infinite, connected, (quasi-)transitive graphs with 
two or more ends have been studied by Li \cite[Thm 1.3]{ZL16}, who
has proved, subject to two conditions, that SAW has positive speed.
The approach of the proof  (see \cite{ZL16}) is to consider a finite `cutset' $W$ with the property that
many SAWs cross $S$
to another component of $G\setminus W$ and never cross back. The pattern theorem is
a key element in the proof.
These results may be applied, for example, to a cylindrical quotient graph of $\ZZ^d$ (see \cite{AJ90,FSG99}),
and the infinite free product of two finite, transitive, connected graphs. Here are two corollaries for Cayley graphs.

\begin{theorem}\cite[Thms 1.7, 1.8]{ZL16}\label{mg}
Let $\Gamma$ be an infinite, finitely generated group with two or more ends.
\begin{letlist}
\item If $\Ga$ has infinitely many ends, and $G$ is a locally finite Cayley graph, there exists $c>0$ such that 
\begin{equation*}
\limsup_{n\to\infty} \bigl| \{\pi\in\Si_n(\id):\|\pi\|\geq c n\}\bigr|^{1/n}=\mu.
\end{equation*}
\item
There exists some locally finite, Cayley graph $G$ such that SAW on $G$ has positive speed.
\end{letlist}
\end{theorem}

\section*{Acknowledgements}
This work was supported in part by the Engineering
and Physical Sciences Research Council under grant EP/I03372X/1.
GRG acknowledges valuable conversations
with Alexander Holroyd concerning Questions \ref{qn:q5} and 
\ref{qn:nonnormal2}, and the hospitality of UC Berkeley during the completion of the work.  
ZL acknowledges support from the
Simons Foundation under grant \#351813, and the National Science Foundation 
under grant DMS-1608896.

\providecommand{\bysame}{\leavevmode\hbox to3em{\hrulefill}\thinspace}
\providecommand{\MR}{\relax\ifhmode\unskip\space\fi MR }
\providecommand{\MRhref}[2]{%
  \href{http://www.ams.org/mathscinet-getitem?mr=#1}{#2}
}
\providecommand{\href}[2]{#2}


\end{document}